\newif\ifreport

\reportfalse
\reporttrue

\newif\ifProofs
\newif\ifStdFormProof

\ifreport

\documentclass[11pt]{article}
\usepackage{array,latexsym,amsmath,amssymb}
\usepackage{amsfonts,cite}
\usepackage{algorithm}
\usepackage{graphicx}
\usepackage{algpseudocode}
\usepackage{caption}
\usepackage{booktabs}
\usepackage{longtable}
\else

\documentclass[onefignum,onetabnum]{siamltex}

\fi

\usepackage{latexsym,amsmath,amssymb}

\DeclareMathOperator*{\Argmax}{Arg\,max}

\date  {February 2014}
\title{AN INEQUALITY-CONSTRAINED SL/QP METHOD FOR MINIMIZING THE SPECTRAL ABSCISSA}
\author{Vyacheslav Kungurtsev\thanks{Agent Technology Center,
        Department of Computer Science,
        Faculty of Electrical Engineering,
        Czech Technical University in Prague,
and Computer Science Department 
and Optimization in Engineering Center (OPTEC), KU Leuven,
Kasteelpark Arenberg 10, B-3001 Leuven- Heverlee, Belgium. 
({\tt vyacheslav.kungurtsev@fel.cvut.cz}).
        }
   \and 
   Wim Michiels \thanks{Computer Science Department and Optimization 
   in Engineering Center (OPTEC), KU Leuven, Kasteelpark Arenberg 10, B-3001 Leuven- Heverlee, Belgium. 
   ({\tt wim.michiels@cs.kuleuven.be} ). } \and
       Moritz Diehl  \thanks{Department of Microsystems Engineering IMTEK,
              University of Freiburg, Georges-Koehler-Allee 102,
              79110 Freiburg, Germany and Electrical Engineering Department and Electrical Engineering Department 
(ESAT-STADIUS) and Optimization in Engineering Center (OPTEC), KU Leuven,
Kasteelpark Arenberg 10, B-3001 Leuven- Heverlee, Belgium. 
({\tt moritz.diehl@imtek.uni-freiburg.de} ). }}

\begin{document}


\Proofsfalse        
\Proofstrue         

\maketitle

\setcounter{page}{1}



\begin{abstract}
We consider a problem in eigenvalue optimization, in particular finding a local minimizer of the 
spectral abscissa - the value of a parameter that results in the smallest value of the largest 
real part of the spectrum of a matrix system. This is an important problem for the stabilization of control 
systems. Many systems require the spectra to lie in the left half plane in order for them to be stable. 
The optimization problem, however, is difficult to solve because the underlying objective function is 
nonconvex, nonsmooth, and non-Lipschitz. In addition, local minima tend to correspond to points of non-differentiability 
and  locally non-Lipschitz behavior. We present a sequential linear and quadratic programming algorithm 
that solves a series 
of linear or quadratic subproblems formed by linearizing the surfaces corresponding to the largest eigenvalues. 
We present numerical results comparing the algorithms to the state of the art.
\end{abstract}




\section{Eigenvalue Problem}\label{sec:problem}


The problem of interest can be written in the form,
\begin{equation}\label{eq:eigprob}
\min_{x\in\Re^n} f(x)=\min_x \alpha(F(x)),
\end{equation}
where the spectral abscissa $\alpha$ is defined to be, 
\[
\alpha(F(x))=\max_i \mathfrak{Re}\left( \lambda_i(F(x))\right),
\]
with $\{\lambda_i(F(x))\}$ is the (possibly infinite) spectrum of the matrix $F(x)$ and 
$F:\Re^n\to \Re^{N\times N}$ is two times continuously differentiable. 
The spectral abscissa corresponds to the largest 
real part of the eigenvalues of $F(x)$. Recall that the set of matrices with 
semi-simple eigenvalues is dense in $\Re^{N\times N}$ and $\frac{d\lambda_i(F(x))}{dx}$ 
is continuous with respect to $x$ for all $i$ and all $x$ such that $F(x)$ has 
only semi-simple eigenpairs, and so $\lambda_i(F(x))$ is locally smooth for a.e. $x$.

For instance, in the field of linear output feedback control, $F(x)$ is defined to be,
\[
F(x) = A+BXC,
\]
where $A$ is the open-loop matrix for the system, $B$ the input matrix and $C$ the
output matrix, and $X$ is formed by arranging the components of $x$ into a matrix of 
the appropriate dimensions. 

The optimization problem is difficult to solve for several reasons: 
\begin{enumerate}
\item It is nonconvex, with possibly many local minimizers and, even arbitrarily close to a 
local minimizer, the spectral abscissa function typically (but not universally) has negative curvature.
\item It is nonsmooth. As the parameter $x$ changes, each eigenvalue changes as well, usually in a 
smooth way, however at points where one smooth eigenvalue surface overtakes another one, points of nonsmoothness 
arise. 
\item It is non-Lipschitz. This occurs in the case of a non semi-simple eigenvalue, at which perturbations of the 
parameter can result in entirely different spectra (for instance, let $x_c$ be such that $F(x_c)$ has a 
non semi-simple eigenvalue, and there exists a complex conjugage eigen-pair for $x>x_c$ and 
two real pairs for $x<x_c$, with the inequalities taken component-wise). Moreover, local minimizers often tend 
to correspond to these points.
\end{enumerate}
These properties induce a challenging task for optimization algorithms. There are many algorithms suitable 
for nonconvex smooth or convex nonsmooth optimization problems, but few that are able to solve 
nonconvex, nonsmooth problems. Furthermore, locally Lipschitz objective functions seem to be universally 
assumed in the analysis of algorithms (such as local convergence of Newton-like methods, 
and sufficient descent and bounded steps of proximal-type first order algorithms). 

On the other hand, the spectral abscissa is a smooth function almost everywhere (a.e.) in any standard measure of 
$\Re^n$, 
so the problem is tractable without necessitating the use of sub-gradients, since gradients can be computed 
at an arbitrary point with probability one. However, local minimizers typically correspond to points of nonsmoothness, 
and so any algorithm seeking a minimizer of $\alpha(F(x))$ should still take into account that the points 
to which the sequence of iterates it generates should be attracted to points at which no unique gradient is defined.

\begin{figure}
\begin{center}
\includegraphics[scale=0.7]{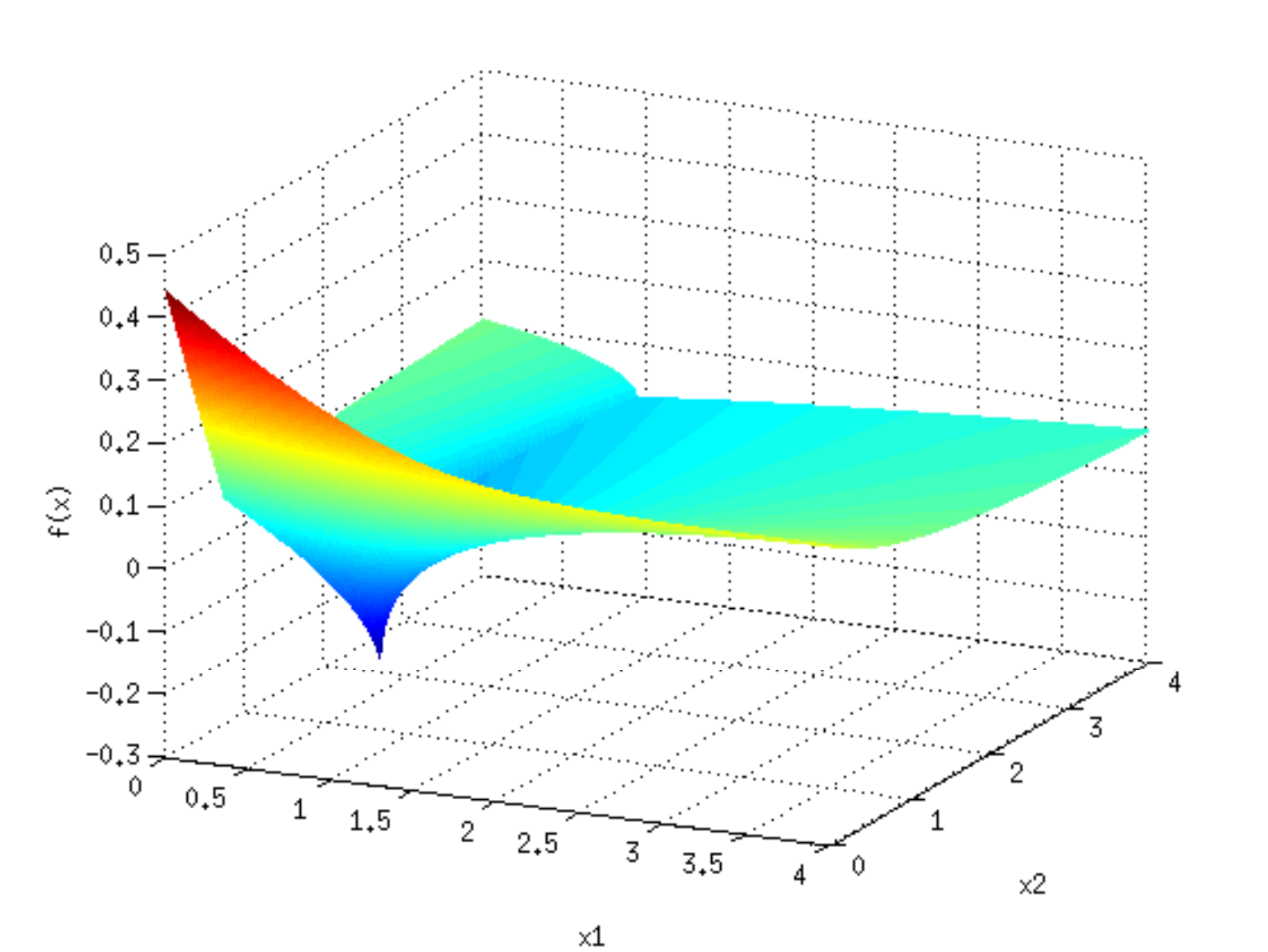}
\end{center}
\caption{\label{fig-plot2d} $f(x)$ for a two-dimensional eigenvalue optimization problem.}
\end{figure}

We present the graph of a two-dimensional problem in Figure~\ref{fig-plot2d}. In this example, 
first given in~\cite{Van2009}, $F(x)=A+BK$, with,
\[
A=\begin{pmatrix} 0.1 & -0.03 & 0.2 \\ 
0.2 & 0.05 & 0.01 \\ 
-0.06 & 0.2 & 0.07 
\end{pmatrix}, \,\,
B=\frac12 \begin{pmatrix} -1 \\ -2 \\ 1 \end{pmatrix}
,\,\, K^T = \begin{pmatrix} x_1 \\ x_2 \\ 1.4 \end{pmatrix} 
\]
Notice that all of 
the features of $\alpha(F(x))$ we describe above, nonconvexity, nonsmoothness and non-Lipchitz behavior
 are evident in the figure. 

\subsection{Previous work}
A number of authors have looked at the problem of achieving matrix stability by minimizing the spectral 
abscissa~\cite{Burke2005, Burke2002, Lewis1996, Burke2006, Apkarian2008, Noll2005, Burke2001}. 
Typically, at each iteration information from a set of gradients is used to generate the next step. 
There are two primary methods of generating and using the appropriate set of gradients. First, in sampling methods, 
a set of gradients is generated by sampling around the current point, 
which serves to approximate the subgradient at a nearby point of nonsmoothness, and a step in the convex hull 
of the sampled gradients is taken~\cite{Burke2005, Burke2002}. 
By contrast, in "bundle" methods, originally developed for convex nonsmooth optimization~\cite{Lemarechal1989} 
build gradient information by mtaintaining historically calculated subgradients during the course of the iterations, 
and at each point solve cutting-plane underapproximations of the function. 
In the case of non-convex problems, the analogue of bundle methods is less 
clear, but various algorithms exist (see, e.g.,~\cite{Fuduli2004}). Finally we mention a unique contribution 
that is not comprehensively described by belonging to either of the two categories, and is specifically 
and algorithm for the context of eigenvalue optimization. 
In~\cite{Noll2005} the $\epsilon$ Clarke-subdifferential is approximated by 
taking the eigenvalues within $\epsilon$ of one with the largest real part and their derivatives with respect to $x$. 
This resembles our approach in using derivative information from multiple eigenvalue surfaces.

There are a limited number of solvers suitable for eigenvalue optimization problems. In this paper we compare 
our results to HANSO~\cite{Overton2009}, a code that uses a combination of a BFGS and a gradient sampling search step. 
BFGS has been shown 
to exhibit good convergence properties for nonsmooth problems~\cite{Lewis2013}. In particular, we have noticed HANSO tends 
to interpret points where the objective function surface locally looks like a cusp (i.e., the function is nonsmooth and 
locally concave in a dense region around the point, the typical situation at a local minimizer of 
$\alpha(F(X))$), the Hessian approximation's largest eigenvalue grows without bound.
The second order method is eventually drawn to such points, which the gradient sampling procedure~\cite{Burke2005} refines. 

\subsection{Contribution}

Recall that the spectral abscissa function $\alpha(F(x))$ is smooth almost everywhere. This is because for each 
eigenvalue, for a.e. $x$, locally there exists a smooth surface $\lambda_i(x)$ respect to $x$ that describes how 
the eigenvalue changes with respect to $x$. For a.e. $x$, $\alpha(F(x))$ corresponds to the surface that is associated 
with the eigenvalue with the maximal real part. Points where $\alpha(F(x))$ is nonsmooth correspond to two situations: 1) a 
switch in which eigenvalue corresponds to the maximal real part as $x$ changes, or 2) a splitting of the eigenvalue 
surface set at a multiple non-semi-simple eigenvalue, in which case the maximal surface can be non-Lipschitz.

 The procedures 
implemented and analyzed thus far only uses the maximal eigenvalue surface, and only after 
extensively sampling regions with other maximal eigenvalues does it refine its steps relative to the 
behavior of the function between these regions. Instead, we consider the possibility of explicitly 
including every one, or at least a large subset, of the eigenvalue surfaces at each point in the calculation of a step 
with a linearized or quadratic approximation. In the case where nonsmoothness corresponds to changes in which 
one eigenvalue surface overtakes another one, this should encourage better steps by anticipating in the 
model where the maximal 
eigenvalue changes, and as a result taking a step that considers both surfaces and, e.g., moving down a 
valley in the topography of $\alpha(F(x))$. In the case of points corresponding to 
a multiple non semi-simple eigenvalue, we introduce memory to include the surface on either side of the 
point of nonsmoothness, which we expect to decrease the number of iterations needed for the optimization 
procedure to take steps that decrease the objective function near these non-Lipschitz points. 

Thus, in constrast to previous methods, which uses sampled or historical linearizations close to the current 
iterate to compute a descent step, the algorithm we present uses linearization information of all, or a selected 
number of the eigenvalue surfaces at the current iterate. We believe that this will allow for larger steps 
to be taken, resulting in a relatively smaller number of iterations required to reach a local minimizer. 

In this paper we include extensive numerical results comparing our algorithm and HANSO. We do not, however, 
include theoretical convergence analysis. We believe that that for the case where the spectral abscissa function 
is locally Lipschitz (which would hold if all eigenvalues of $F(x)$ are semi-simple for all $x$), then such a 
proof would be simple and straightforward by standard arguments. However, in the general case, minimization of 
non-locally-Lipschitz functions is a very difficult problem with respect to convergence theory. We note that 
the convergence proof for gradient sampling~\cite{Burke2005} relies on the local Lipschitz property.

\subsection{Notation}
We will assume that the eigenvalues are ordered by their real components, with ties broken arbitrarily, i.e., 
\[
\mathfrak{Re}(\lambda_0(F(x))) \ge \mathfrak{Re}(\lambda_1(F(x))) \ge \mathfrak{Re}(\lambda_2(F(x))) \ge ...  
\]
note that this ordering is dependent on $x$.

\section{SL/QP for Eigenvalue Optimization}\label{sec:algorithm}
In this section we present a sequential quadratic programming (SQP) method, as well as a more 
simplified Sequential Linear Programming (SLP) method for minimizing the spectral absicssa. 
The algorithm is a trust-region based method based on the realization 
that the problem $\min_{x\in\Re^n} \alpha(F(x))$ can be rewritten as,
\begin{equation}\label{eq:eigprobnew}
\begin{array}{rl}
\min_{\gamma\in\Re,x\in\Re^n} & \gamma, \\
\text{subject to} & \gamma \ge \mathfrak{Re}(\lambda_i(F(x))) \text{ for all } i.
\end{array}
\end{equation}

One key observation to inspire the algorithm is that the number of points at which 
the objective function $\alpha(F(x))$ in~\eqref{eq:eigprob} is nonsmooth is of measure zero in the Lebesgue space 
$\Re^n$. This implies that for a.e. $x$, the function $\alpha(F(x))$ is a locally smooth surface. This surface 
corresponds to the value of $\lambda_0(F(x))$ as a function of $x$. 

Locally, we can express a linear approximation of the spectral abscissa function as a plane 
through the point $(x,\alpha(F(x)))$ with the gradient $\nabla_x \alpha(F(x))=
\nabla_x\lambda_0(F(x))$, if $\lambda_0(F(x))$ is simple. 
We can calculate both this vector as well as $\nabla^2_{xx} \lambda_i (F(x))$ by the formulas~\cite{Lancaster1964,Magnus1985}, 
\begin{equation}
\nabla_x \lambda_i (F(x)) = \frac{u_i^* \frac{dF}{dx} v_i}{u_i^* v_i},
\end{equation}
\begin{equation}
\label{eq:secondderiv}
\nabla^2_{xx} \lambda_i (F(x)) = \frac{2 u_i^* \left(\frac{d F}{dx}\right) K_i(\lambda_i I-F(x))^\dagger K_i 
\left(\frac{d F}{dx}\right) v_i}{u^*_i v_i},
\end{equation}

where $u_i$ and $v_i$ are the left and right eigenvectors of $F(x)$ corresponding to eigenvalue $i$, $u^*$ 
corresponds to the conjugate of $u$, $\dagger$ is the pseudo-inverse, and $K=I-v_i u_i^*/(u^*_i v_i)$. 

The Lagrangian function for the problem~\eqref{eq:eigprobnew} is defined as,
\begin{equation}
\label{eq:lag}
L(\gamma,x) = \gamma-\sum_i y_i (\gamma- \mathfrak{Re}\left(\lambda_i (F(x))\right)),
\end{equation}
where $y$ is the vector of Lagrange multipliers.

This naturally suggests the SQP method wherein a sequence of iterations $x_{k+1}=x_k+t\Delta x$ is calculated, 
with $t$ a line-search scalar and $\Delta x$ is determined by solving subproblems of the form,
\begin{equation}\label{eq:basicsqp}
\begin{array}{rl}
\min_{\Delta x,\Delta\gamma} & \Delta \gamma +\frac12\Delta x^T H_k 
\Delta x, \\
\text{subject to} & \Delta\gamma+\alpha(F(x_k)) 
\ge \mathfrak{Re}(\lambda_i(F(x_k)))+\mathfrak{Re}(\nabla_x\lambda_i(F(x_k)))^T \Delta x, \forall i,
\end{array}
\end{equation}
for $\Delta x$ and $\Delta \gamma$, where $H_k$ is a Lagrangian Hessian term at $x_k$. 
This resembles the minimax subproblem arising in bundle methods~\cite{Lemarechal1989}, except 
that the linearizations are defined around a set of eigenvalues of $F(x)$ evaluated at one $x_k$.
At each iteration we discard $\Delta \gamma$ and compute $\gamma_{k+1} = \alpha(F(x_{k+1}))$ instead and so the 
iteration procedure behaves like a slack reset in nonlinear programming. 
To determine the Lagrangian Hessian, we do not use the multipliers from the quadratic programs 
as estimates for the Lagrangian multipliers. This is because the linearizations are local, 
and the surfaces they correspond to may disappear from one point to the next (in case of a non-semi-simple 
eigenvalue between them), and so each constraint may not correspond in any meaningful way to a constraint 
estimated previously. 
In the sense of problem~\eqref{eq:eigprobnew}, the eigenvalues in the set $\Argmax_{i=0,1,...} (\mathfrak{Re}(\lambda_i(F(x))))$ 
are the \emph{active} constraints. For a.e. $x$ this set comprises of at most two elements, corresponding to a 
conjugate pair, but in this case the surface corresponding to the real value of the eigenvalue as a function of $x$ is
the same for each eigenvalue of the pair.
So at the start of each iteration, we set the Lagrange multiplier to have $1$ in the component corresponding to  
the/an eigenvalue with maximal real part, and $0$ otherwise. Using this multiplier, we calculate the Hessian,
\[
H_k = \mathfrak{Re}\nabla^2_{xx} \lambda_0(F(x_k))).
\]
Note that the Lagrangian 
function~\eqref{eq:lag} is linear with respect to $\gamma$ and so the Hessian only has blocks corresponding 
to $x$. At each step of the SQP we need only update $x_{k+1}=x_k+\Delta x$. The solution $\Delta \gamma$ is discarded.

There are two primary additional features to the basic procedure we have presented in order to make the method more 
practically successful. First, recall from Section~\ref{sec:problem}, that the underlying problem is non-convex. This implies 
that at any local quadratic approximation of an eigenvalue surface, the Hessian could be indefinite or even negative 
definite. This implies that the approximating quadratic program~\eqref{eq:basicsqp} could be unbounded below. We 
constrain the problem with a trust-region to prevent this. Since we have linear constraints, we use an infinity norm trust-region, 
which acts as a "box" limiting the magnitude of the maximal component of $\Delta x$. 
\begin{equation}\label{eq:basicsqpandtr}
\begin{array}{rl}
\min_{\Delta x,\Delta\gamma} & \Delta \gamma +\frac12 \Delta x^T H_k 
 \Delta x  , \\
\text{subject to} & \Delta\gamma+\alpha(F(x_k)) \ge \mathfrak{Re}(\lambda_i(F(x_k)))+\mathfrak{Re}(\nabla_x\lambda_i(F(x_k)))^T \Delta x, \, \forall i, \\
& ||\Delta x||_\infty \le \Delta_k.
\end{array}
\end{equation}
Since the Hessian could be indefinite, the solution $\Delta x$ could be a direction of ascent for the objective function. 
If the second order information is sufficiently accurate then the step should still decrease the objective function. 
Otherwise, however, it could be that any point along the line segment $x_k+t\Delta x$, $t\in [0,1]$ satisfies
$\alpha(F(x_k+t\Delta x))>\alpha(F(x_k))$. Hence, after computing $\Delta x$ we first test if,
\begin{equation}\label{eq:acceptfull}
\alpha(F(x_k+\Delta x))<\alpha(F(x_k)),
\end{equation}
in which case we set $x_{k+1}=x_k+\Delta x$ and continue to the next iteration. Otherwise, we test for descent, 
\begin{equation}\label{eq:testdescent}
\mathfrak{Re}(\nabla_x\lambda_i(F(x_k)))^T \Delta x < 0,
\end{equation}
and if this does not hold we set $\Delta_{k+1} = \gamma_1 \Delta_k$, where $\gamma_1$ is a constant satisfying $\gamma_1\in(0,1)$, 
and resolve the subproblem.

If~\eqref{eq:testdescent} holds, we follow the mixed trust-region/line-search procedure presented by Gertz~\cite{Gertz1999}, 
in which a backtracking line search reduces the size of the 
step $t$ until decrease is achieved $(\alpha(F(x_k+t \Delta x))<\alpha(F(x_k)))$, and the next trust-region radius corresponds to $t ||\Delta x||$. 

\begin{equation}\label{eq:updatedelta}
\Delta_{k+1}=\left\{
\begin{array}{lr}
\gamma_2 \Delta_k & \text{ if } \alpha(F((x_k+\Delta x)) < \alpha(F(x_k)) \\
t ||\Delta x|| &  \text{ otherwise},
\end{array}\right.
\end{equation} 
where $\gamma_2$ is a constant satisfying $\gamma>1$.

We update the trust-region simply by increasing it if we achieve descent, and decreasing it otherwise. For consistency 
with convergence theory~\cite{Conn2000}, we would enforce sufficient decrease conditions with respect to predicted (from the 
quadratic approximation) and actual decrease. However, since lax criteria of acceptance (e.g., with a small constant multiplying 
the predicted-actual decrease ratio)
of the step is practically equivalent to this condition, we proceed as in the line-search criteria for the gradient 
sampling method~\cite{Burke2005} to just enforce descent.

\begin{figure}
\begin{center}
\includegraphics[scale=0.7]{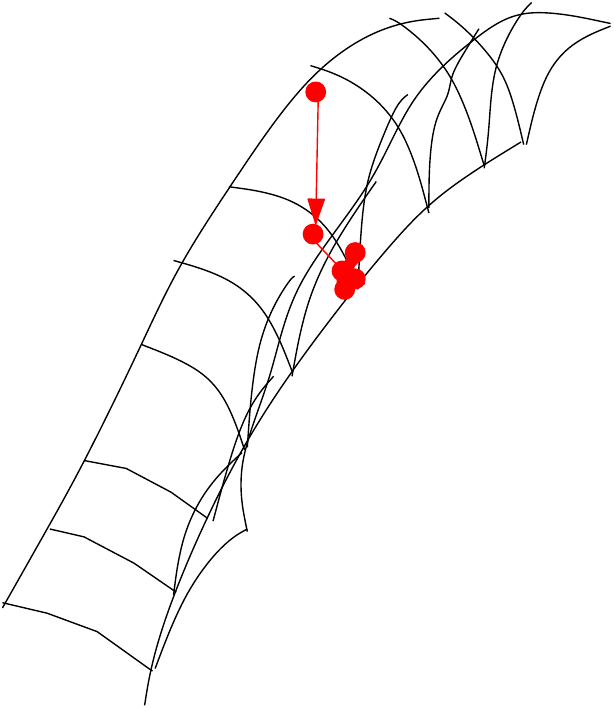}
\end{center}
\caption{\label{fig-valley} Possible set of iterations of SQP without memory along a surface of $\alpha(F(x))$.}
\end{figure}

In addition, recall that near a non semi-simple eigenvalue, 
there could be qualitatively different eigenvalue surface combinations on either side of a "valley", or $n-1$ dimensional 
hypersurface in $\Re^n$ at which $\nabla \lambda_i(F(x))$ is undefined. 
In order to account for this, we added another feature to the algorithm, after observing a certain phenomenon 
that was typical with the original SQP algorithm with the trust region but without this additional feature. 
In many cases, the algorithm jammed near valleys of this kind and would frequently converge 
onto the valley rather than move down along it. This is because locally, the directional derivative of $\alpha(F(x))$ is steeper  
towards the valley than perpendicular to it, and, furthermore, it has negative curvature along that direction, 
so a local approximation that regards only the eigenvalue surfaces at a point on one side of the valley will 
result in the step of steepest decrease being in this direction. Since the surface on the other side of the valley 
is not accounted for on the original side, this is not incorporated directly into the subproblem. We illustrate 
this scenario in Figure~\ref{fig-valley}.

To remedy this, we added "memory" to the SQP method with a set $\mathcal{M}$. When it occurs that 
$\alpha(F(x_k+\Delta x)) > \alpha(F(x_k))$, which we expect in the "jamming" scenario described above, 
the procedure stores the tuple $\{ x_k+\Delta x_k,\mathfrak{Re}(\lambda_0(F(x_k+\Delta x_k))),
\nabla_x \mathfrak{Re}(\lambda_0(F(x_k+\Delta x_k)))\}$ in $\mathcal{M}$.
Then, if in a future iteration $l$, the current point $x_l$ satisfies $||x_l-x^{i}||_\infty \le \Delta_k$ for any $i\in\{1,...,|\mathcal{M}|\}$, 
then we include this linearized surface in the QP subproblem. 

\begin{equation}\label{eq:basicsqpandtrm}
\begin{array}{rl}
\min_{\Delta x,\Delta\gamma} & \Delta \gamma +\frac12 \Delta x^T H_k \Delta x, \\
\text{subject to} & \Delta\gamma+\alpha(F(x_k)) \ge \mathfrak{Re}(\lambda_i(F(x_k)))+\mathfrak{Re}(\nabla_x\lambda_i(F(x_k)))^T \Delta x, \, \forall i \\
& \Delta\gamma+\alpha(F(x_k)) \ge \mathfrak{Re}(\lambda_{(i)}(F(x^{(i)})))
\\ & \qquad \qquad \qquad \qquad +\mathfrak{Re}(\nabla_x\lambda_{(i)}(F(x^{(i)})))^T (x_k+\Delta x-x^{(i)}), \, i\in M_k \\ 
& ||\Delta x||_\infty \le \Delta_k,
\end{array}
\end{equation}
where $M_k\subset \mathcal{M}$ represents the points $x^{(i)}$ satisfying $||x^{(i)}-x_k||\le \Delta_k$.

It can be pointed out that this resembles gradient bundle algorithms, which uses gradient information 
from previous iterations to generate the current step. We acknowledge the resemblance, but point out that 
in this case historical gradient information is only added 
in the limited case of a failed trial step, rather than used in each iteration, it is only used when 
the underlying eigenvalue surfaces change. Local eigenvalue surface approximations, however, are used
at each iteration. 

Finally, we consider two simplifications of the subproblem. First, we consider dropping the Hessian term, to 
formulate a sequential linear programming procedure. We note that the lack of a Lipschitz property makes local 
quadratic convergence difficult to prove for any SQP/Newton type scheme, and so it is quite 
possible that second derivatives will not reduce the computation 
time. In addition, with possibly nonconvex subproblems, QP subproblems could result in a series of 
ascent steps requiring resolving with a smaller trust-region, whereas with LP subproblems, descent 
is guaranteed for any solution.

We present the SLP version below~\eqref{eq:basicslpandtr} and summarize the full algorithm as 
Algorithm~\ref{alg-slpeval}, for simplicity, 
and discuss a comparison of the SLP and SQP variant in the numerical results 
section. Note that we can select a subset of the eigenvalues $N_k$ at each iteration $k$ 
to evaluate and linearize. For small problems, we set $N_k=N$ for all $k$ (where, recall that
$F(x)\in\Re^{N\times N}$).

\begin{equation}\label{eq:basicslpandtr}
\begin{array}{rl}
\min_{\Delta x,\Delta\gamma} & \Delta \gamma , \\
\text{subject to} & \Delta\gamma+\alpha(F(x_k)) \ge \mathfrak{Re}(\lambda_i(F(x_k)))
\\ & \qquad \qquad \qquad \qquad+\mathfrak{Re}(\nabla_x\lambda_i(F(x_k)))^T \Delta x, \, i\in\{0,...,N_k\} \\
& \Delta\gamma+\alpha(F(x_k)) \ge \mathfrak{Re}(\lambda_{(i)}(F(x^{(i)})))
\\ & \qquad \qquad \qquad \qquad+\mathfrak{Re}(\nabla_x\lambda_{(i)}(F(x^{(i)})))^T (x_k+\Delta x-x^{(i)}), \, i\in M_k \\ 
& ||\Delta x||_\infty \le \Delta_k.
\end{array}
\end{equation}

The stopping criterion corresponds to the step becoming small, without any new information (memory) being added at 
the current iteration.

\begin{algorithm}[H]
\caption{\label{alg-slpeval} SLP Algorithm for Eigenvalue Optimization.}
\begin{algorithmic}[1]
\State Define constants $0<\gamma_1<1$, $\gamma_2>1$, $\delta_m>0$, and $S\in \mathbb{N}$. 
\State Determine $N_0$. Typically, set $N_0=N$, the size of $F(x)$.
\For{$S$ times}
\State Randomly select starting point $x_0$.
\State Set $\mathcal{M}_1=\emptyset$. Set $k=1$. 
\State Calculate initial $\{\lambda_i(F(x_0))\}$ and $\{\nabla_x \lambda_i(F(x_0))\}$ for $i\in\{0,..,N_0\}$.
\While{($||\Delta x|| > \delta_m$ or $\mathcal{M}_k\neq\mathcal{M}_{k-1}$)}
\State Solve~\eqref{eq:basicslpandtr} for $\Delta x_k$.
\State Calculate $\{\lambda_i(F(x_k+\Delta x_k))\}$ and $\{\nabla_x \lambda_i(F(x_k+\Delta x_k))\}$ for $i\in\{0,..,N_k\}$
\If{$\alpha(F(x_k+\Delta x_k))<\alpha(F(x_k))$}
\State Set $x_{k+1}\gets x_k$
\State Set $\Delta_{k+1} \gets \gamma \Delta_k$.
\Else
\State Store $\{x_k+\Delta x_k,\mathfrak{Re}(\lambda_0(F(x_k+\Delta x_k))),
\nabla_x \mathfrak{Re}(\lambda_0(F(x_k+\Delta x_k)))\}$ 
\State \qquad \qquad \qquad in $\mathcal{M}_{k+1}$.
\State Find $t$ such that $\alpha(F(x_k+t \Delta x_k))<  \alpha(F(x_k))$.
\State Set $x_{k+1} \gets x_k+t \Delta x_k$. 
\State Set $\Delta_{k+1} \gets t ||\Delta x_k||$.
\EndIf
\State Set $k\gets k+1$.
\State Determine $N_k$. Typically, set $N_k=N$, the size of $F(x)$.
\State Calculate all $\{\lambda_i(F(x_k))\}$ and $\{\nabla_x \lambda_i(F(x_k))\}$ for $i\in\{0,..,N_0\}$.
\EndWhile
\State Add the last point $(x_f,\alpha(F(x_f)))$ to $\mathcal{F}$.
\EndFor 

\Return{$\{x_f,\alpha(F(x_f))\}$ corresponding to the lowest value of $\alpha(F(x_f))$ in $\mathcal{F}$.}
\end{algorithmic}
\end{algorithm}

In the case of large scale problems, we reduce $N_k$ in order to minimize eigenvalue computation time. 
For large $N$ we use iterative procedures to compute a subset of the eigenvalues corresponding to those 
farthest to the right in the complex plane. For example, want to calculate all eigenvalues in the right half 
plane relative to a certain minimal value of some desired $\bar{\lambda}_c$. 
The set $N_k$ of eigenvalues we consider may change each iteration. 
We can choose to either determine $\bar{\lambda}_c$ or $N_k$ and calculate the appropriate subset of 
the spectrum. If we decide to select the subset based on a heuristic choice for $N_k$, its 
value should be greater than or equal to the number of active (maximal) eigenvalues, otherwise the
choice depends on how many surfaces we want to keep track of, locally. In our implementation we use $2n$.
Note that this is the number of sample points 
to use in a gradient sampling method~\cite{Burke2005}

\section{Numerical Results for Linear Eigenvalue Problems}\label{s:complib}
We compare SLP/SQP to HANSO on a set of linear control problems arising from COMPlib~\cite{Leibfritz2010}. 
COMPlib contains a set of linear 
control matrices, of which we take three for the system $F(x) = A+BXC$. Out of 124 examples, we picked 99 problems with 
the number of rows (and columns) of $A$ was less than or equal to 50. We found that for both HANSO and SLP/SQP, even 
for the large-scale variant, many of the problems in the complement of this set took far long to solve to 
realistically perform numerical comparisons.

For all solvers, we used a stopping tolerance of 1e-6 (indicating that the algorithms stopped when the (inf) norm of 
the step was smaller than 1e-4). The SL/QP algorithms were coded in MATLAB, with all 
tests run using MATLAB version 2013a. All tests were performed on an Intel Core 2.2 GHz $\times 8$ running Ubuntu 14.04.
For all algorithms we use the same procedure of using ten random starting points as provided by default 
with HANSO, specifically initializing a point by a normal distribution centered at zero, and then picking the 
best solution (the one with the lowest objective value) of ten runs. All of the code used for the experiments in 
this and the next section is available at~\cite{SLPCode}.

We list the parameter and initial values we use in our implementations of SL/QP in Table~\ref{tab:params}
$0<\gamma_1<1$, $\gamma_2>1$, $\Delta_m>0$, $\delta_m>0$, and $S\in \mathbb{N}$. We denote $k_{\text{max}}$ 
the maximum number of iterations, $LSk_{\text{max}}$ the maximum number of line-search steps, and 
$\eta$ the backtracking contraction parameter.
\begin{table}[h]
\caption{Control parameters and initial values required by algorithm~\ref{alg-slpeval}}\label{tab:params}
\begin{footnotesize}
\begin{center}
  \begin{tabular}{|>{\tt}c>{\tt}l|>{\tt}c>{\tt}l|>{\tt}c>{\tt}l|} \hline
\multicolumn{1}{|c}{\rm Parameter} &
\multicolumn{1}{c}{\rm Value} &
\multicolumn{1}{|c}{\rm Parameter} &
\multicolumn{1}{c}{\rm Value} &
\multicolumn{1}{|c}{\rm Parameter} &
\multicolumn{1}{c|}{\rm Value} \\ \hline
 $\gamma_1$           & 0.1    & $S$                        & 10 & $k_{\text{max}}$ & 20\\
 $\gamma_2$    & 2.0 & $N_0$                         & $N$    & $LSk_{\text{max}}$      & 20 \\
 $\delta_m$        & 1.0e-r    & $\Delta_0$                      & 1.0 & $\eta$     & 0.5    \\\hline
\end{tabular}
\end{center}
\end{footnotesize}
\end{table}

HANSO was used with its default parameters, including {\tt normtol}$=$1e-4, {\tt evaldist}$=$1e-4, and 
{\tt maxit}$=$1000.

For these comparisons, we compare both the time of execution as well as the value of the final solution. 
We summarize the results below in Table~\ref{tablecomplib}. If an algorithm returns an error rather than 
converging, we indicate that as being the worse performer (for no problem did both SLP and a variation of HANSO 
fail to converge).
Note that gradient sampling provides a guarantee of 
global convergence in the case of the objective being local Lipschitz, however there are no theoretical guarantees 
for any algorithm without this condition. 

Since the problems are all nonconvex, it is difficult to make a straightforward comparison since each run could result 
in a different local minimizer, and thus we perform a large number of runs for each problem to obtain a global 
picture. Since for all of the algorithms we sample ten random starting points from a uniformly normal distribution, 
in the long run, since the algorithms are all purely local in nature (i.e., encouraging convergence to a local minimizer 
from arbitrary starting points, rather than a global minimizer), a global picture of which algorithm performs 
more reliably in obtaining a local minimizer can be understood as one for which a) convergence from arbitrary starting 
points is more likely and b) possibly includes general features that encourages a globally lower objective value.
In the results we report the median, minimum, and maximum number of problems for which SLP outperforms HANSO from a set of 50 
runs of every small problem.

In general, for most problems, SLP 
finds a lower minimum than HANSO. It appears to be faster than HANSO with gradient sampling, and slower than 
BFGS alone. Interestingly, gradient sampling does not, on average, tend to improve the performance of HANSO 
vis-a-vis SLP.

In the interest of enforcing stability, for SLP, for one run there were 74 problems 
where the final spectral abscissa (median of 50 runs) was less than zero, as compared to 63 for HANSO with 
gradient bundle and 59 for HANSO without. 

\begin{table}
\caption{\label{tablecomplib} Number of times SLP outperformed HANSO out of 99 small problems. 
Minimum (median) and maximum out of 50 runs. The value for each algorithm is the best out of 10 
random starting points for each run and each problem. The 
time is the total clock time taken to run the algorithm for all of these random starting points.}
\begin{center}
\begin{tabular}{l  c  r } 
\toprule 
 & in value  & in time \\ \midrule 
HANSO & 68 (78) 85 & 85 (90) 93   \\ 
HANSO without gradient sampling & 67 (78) 85 & 4 (8) 12  \\

\bottomrule
\end{tabular}
\end{center}
\end{table}

We also performed a comparison of SLP versus SQP on the test problems. The results, in 
Table~\ref{tablecomplibslqp} indicate that SLP tends to be slightly more reliable and faster. 
Overall, there is little difference in the performance, however.
\begin{table}
\caption{\label{tablecomplibslqp} Performance of SLP as compared to SQP (out of 99 problems, minimum (median) and maximum
 of 50 runs). The value for each algorithm is the best out of 10 random starting points for each run and each problem. The 
time is the total clock time taken to run the algorithm for all of these random starting points.
(Notice that the sum of the best in value column is not (and is not very close to) 
$99$, this is because for some problems the result was exactly equal)}
\begin{center}
\begin{tabular}{l  c  r } 
\toprule 
Algorithm & best in value  & best in time \\ \midrule 
SLP & 41 (51) 61 & 49 (58) 66   \\ 
SQP & 35 (46) 54 & 33 (42) 50  \\
\bottomrule
\end{tabular}
\end{center}
\end{table}

Table~\ref{tab.bigtable} lists all of the problems, times and values for SLP and the two variations of HANSO. 
NaN indicates that the solver did not converge for the problem, with either non-convergence or an error indicated. 
We can see that when there are differences in final value, that the differences are typically relatively substantial. 
Indeed, note that, in general, it clearly appears that differences in value must correspond to either different local 
minimizers or failure to converge, rather than a slightly more precise solution.

\begin{longtable}{|>{\tt}l
        @{\hspace{10pt}}>{\tt}r@{\hspace{10pt}}>{\tt}r@{\hspace{10pt}}>{\tt}r@{\hspace{10pt}}>{\tt}r
        @{\hspace{10pt}}>{\tt}r@{\hspace{10pt}}>{\tt}r@{\hspace{10pt}}>{\tt}r@{\hspace{10pt}}>{\tt}r|}
      \caption{\label{tab.bigtable} Results for SLP, HANSO, and HANSO without gradient sampling for 99 small COMPleib problems. Value 
is the best out of 10 runs with random starting points. Time is the total time taken to perform the 10 runs.}\\\hline
\multicolumn{1}{|l}{\rm Prob}&$N$&$n$&
\multicolumn{1}{c}{\rm t SLP}&
\multicolumn{1}{c}{\rm t HANSO}&
\multicolumn{1}{c}{\rm t H no GS}&
\multicolumn{1}{c}{\rm Value SLP}&
\multicolumn{1}{c}{\rm V HANSO}&
\multicolumn{1}{c|}{\rm V H no GS}\\\hline
\endhead
\hline
\endfoot
1 & 5 &  9 & 4.85e-01 &  8.50e-01  & 1.35e+01 &  -1.33e-01 &  -4.99e-01  & -2.05e-01 \\ 
2 & 5 &  9 & 6.83e-01 &  8.44e-01  & 1.33e+01 &  -3.23e-01 &  -4.14e-01  & -8.81e-02 \\ 
3 & 5 &  8 & 1.08e+00 &  8.89e-01  & 6.96e-01 &  -1.39e+00 &  -9.93e-01  & -1.35e-01 \\ 
4 & 4 &  2 & 2.28e-01 &  3.40e-01  & 6.73e-02 &  -5.00e-02 &  -5.00e-02  & -5.00e-02 \\ 
5 & 4 &  4 & 4.36e-01 &  6.30e-01  & 6.37e+00 &  -8.99e-01 &  9.66e-01  & 9.94e-01 \\ 
6 & 7 &  8 & 7.80e-01 &  5.20e-01  & 1.34e+01 &  -6.46e-01 &  -7.53e-01  & -2.17e-01 \\ 
7 & 9 &  2 & 3.53e-01 &  3.07e-01  & 2.16e-01 &  -3.62e-02 &  -3.34e-02  & -3.09e-02 \\ 
8 & 9 &  5 & 7.09e-01 &  1.08e+00  & 8.41e+00 &  -2.50e-01 &  -3.74e-01  & 3.36e-02 \\ 
9 & 10 &  20 & 6.45e-01 &  9.61e-01  & 1.41e+01 &  -2.96e-01 &  -1.09e-01  & -3.04e-02 \\ 
11 & 5 &  8 & 1.17e+00 &  1.33e+00  & 9.39e+00 &  -8.76e+00 &  -2.57e+00  & -1.17e+00 \\ 
12 & 4 &  12 & 6.10e-01 &  8.98e-01  & 1.23e+01 &  -2.63e-01 &  -1.07e-01  & -9.38e-02 \\ 
13 & 28 &  12 & 5.70e-01 &  9.11e-01  & 7.00e-01 &  1.44e-01 &  1.33e-01  & 1.42e-01 \\ 
14 & 40 &  12 & 8.74e-01 &  9.45e-01  & 2.24e+01 &  3.32e-02 &  1.92e-01  & 2.22e-01 \\ 
15 & 4 &  6 & 5.24e-01 &  9.77e-01  & 7.56e+00 &  -4.91e-01 &  -3.85e-01  & -1.46e-01 \\ 
16 & 4 &  8 & 1.17e+00 &  1.10e+00  & 9.50e+00 &  -7.88e-01 &  -4.84e-01  & -1.68e-01 \\ 
17 & 4 &  2 & 3.59e-01 &  5.90e-01  & 5.42e+00 &  -9.35e-01 &  -1.13e+00  & -8.92e-01 \\ 
18 & 10 &  4 & NaN &  3.33e-01  & 7.89e+00 &  NaN &  1.89e+00  & 8.89e-01 \\ 
19 & 4 &  2 & 4.70e-01 &  2.85e-01  & 5.48e+00 &  -2.39e-01 &  -6.27e-02  & -1.49e-01 \\ 
20 & 4 &  4 & 3.86e-01 &  3.60e-01  & 6.52e+00 &  -7.03e-01 &  -1.66e-01  & -6.26e-01 \\ 
21 & 8 &  24 & 5.70e-01 &  1.02e+00  & 1.88e+01 &  -2.87e-02 &  -6.76e-02  & 1.49e-01 \\ 
22 & 8 &  24 & 1.72e+00 &  1.46e+00  & 1.93e+01 &  -2.43e-01 &  -1.08e-01  & 5.47e-02 \\ 
23 & 8 &  8 & 1.09e+00 &  8.30e-01  & 9.80e+00 &  -1.60e-02 &  -6.33e-02  & 1.76e-01 \\ 
24 & 20 &  24 & 5.34e-01 &  1.76e+00  & 8.18e-01 &  -5.00e-03 &  -5.00e-03  & 3.93e-02 \\ 
25 & 20 &  24 & 7.49e-01 &  1.37e+00  & 4.65e+00 &  -5.00e-03 &  -5.00e-03  & 2.04e-02 \\ 
26 & 30 &  15 & 1.18e+00 &  1.41e+00  & 4.69e-01 &  5.53e+00 &  9.36e+00  & 9.79e+00 \\ 
27 & 21 &  9 & 6.06e-01 &  8.94e-01  & 2.02e+01 &  -3.82e-01 &  -4.71e-01  & 1.98e-01 \\ 
28 & 24 &  18 & 3.99e-01 &  2.68e+00  & 2.78e+01 &  -6.22e-01 &  -1.59e+00  & 2.36e+00 \\ 
29 & 4 &  6 & 5.62e-01 &  1.20e+00  & 9.52e+00 &  6.89e-02 &  -2.57e+00  & -1.84e+00 \\ 
30 & 4 &  4 & 5.05e-01 &  4.25e-01  & 6.84e+00 &  -2.11e+00 &  -1.03e+00  & -1.11e+00 \\ 
31 & 12 &  3 & 7.48e-01 &  4.52e-01  & 7.14e+00 &  -2.07e-02 &  -2.07e-02  & -2.07e-02 \\ 
32 & 8 &  1 & 2.27e-01 &  2.28e-01  & 5.74e+00 &  7.12e-01 &  6.63e-01  & 8.66e-01 \\ 
33 & 8 &  16 & 8.75e-01 &  1.22e+00  & 2.02e+01 &  -7.99e-01 &  -7.75e-01  & -3.93e-01 \\ 
34 & 3 &  4 & 1.67e+00 &  9.95e-01  & 6.31e+00 &  -8.52e+00 &  -1.29e+01  & 4.49e-01 \\ 
35 & 6 &  16 & 1.63e+00 &  1.17e+00  & 1.47e+01 &  -4.48e+00 &  -2.36e+00  & -3.81e-01 \\ 
36 & 6 &  24 & 1.00e+00 &  1.37e+00  & 1.87e+01 &  -1.85e+01 &  -2.21e+00  & -1.02e+00 \\ 
37 & 4 &  4 & 7.60e-01 &  9.75e-01  & 6.57e+00 &  -2.16e+00 &  -1.29e+00  & 8.93e-01 \\ 
38 & 10 &  4 & 6.79e-01 &  5.70e-01  & 7.75e+00 &  -7.06e-01 &  -7.11e-01  & 1.17e-01 \\ 
39 & 12 &  4 & 3.02e-01 &  5.91e-01  & 8.06e+00 &  -2.16e-01 &  -2.16e-01  & -2.16e-01 \\ 
40 & 10 &  12 & 4.15e-01 &  4.15e-01  & 1.47e+01 &  -1.71e-01 &  -4.36e-01  & -5.23e-01 \\ 
41 & 10 &  12 & 5.06e-01 &  6.07e-01  & 1.76e+01 &  -1.56e+00 &  -1.47e+00  & -6.70e-01 \\ 
42 & 10 &  12 & 4.79e-01 &  5.52e-01  & 1.96e+01 &  -1.10e+00 &  -4.29e+00  & -7.86e-01 \\ 
44 & 11 &  9 & 9.85e-01 &  1.36e+00  & 7.19e-01 &  -3.62e-03 &  -6.58e-03  & -3.22e-03 \\ 
46 & 4 &  6 & 8.57e-01 &  6.83e-01  & 8.08e+00 &  -2.65e-02 &  -3.13e-02  & -2.34e-02 \\ 
47 & 8 &  4 & 7.63e-01 &  1.11e+00  & 7.40e+00 &  1.78e+01 &  -1.49e+01  & -1.17e+01 \\ 
48 & 21 &  110 & 9.76e-01 &  1.94e+01  & 7.22e-01 &  1.03e+00 &  1.13e+00  & 1.22e+00 \\ 
49 & 20 &  20 & 6.65e-01 &  1.17e+00  & 2.94e+01 &  -9.28e-02 &  -1.03e-01  & -1.12e-01 \\ 
51 & 10 &  1 & 2.07e-01 &  2.79e-01  & 6.22e+00 &  -1.39e-01 &  -1.15e-01  & -1.23e-01 \\ 
52 & 10 &  1 & 3.15e-01 &  2.86e-01  & 6.16e+00 &  -1.48e-01 &  -1.48e-01  & -1.13e-01 \\ 
53 & 10 &  1 & 2.22e-01 &  2.65e-01  & 6.22e+00 &  -9.56e-02 &  -1.09e-01  & -9.84e-02 \\ 
54 & 20 &  1 & 1.45e+00 &  5.37e-01  & 1.10e+01 &  -8.52e-02 &  -8.52e-02  & -2.02e-02 \\ 
55 & 40 &  1 & 3.15e+00 &  1.65e+00  & 9.88e-01 &  -4.02e-05 &  -4.02e-05  & -4.01e-05 \\ 
57 & 5 &  3 & 1.85e-01 &  1.15e+00  & 5.95e+00 &  -5.72e-07 &  -1.33e-06  & -7.57e-07 \\ 
58 & 7 &  8 & 1.30e+00 &  1.21e+00  & 9.94e+00 &  -6.83e-02 &  -2.25e-02  & 9.38e-03 \\ 
59 & 7 &  6 & 4.32e-01 &  6.10e-01  & NaN &  -1.00e-05 &  -1.00e-05  & NaN \\ 
60 & 7 &  6 & 6.10e-01 &  5.95e-01  & 7.90e+00 &  3.88e-02 &  2.91e-03  & 1.21e-01 \\ 
61 & 7 &  6 & 5.91e-01 &  3.67e-01  & 8.37e+00 &  -2.11e+00 &  -1.90e+00  & -1.21e+00 \\ 
64 & 3 &  2 & 2.52e-01 &  4.42e-01  & 1.32e+00 &  1.45e+00 &  1.48e+00  & 2.84e+00 \\ 
65 & 2 &  1 & 5.19e-01 &  6.89e-01  & 4.75e+00 &  -1.00e+00 &  -1.00e+00  & -9.60e-01 \\ 
66 & 4 &  1 & NaN &  3.47e-01  & 2.17e-01 &  NaN &  2.18e+00  & 2.14e+00 \\ 
67 & 4 &  6 & 7.93e-01 &  3.87e-01  & 7.50e+00 &  -1.31e+00 &  -6.17e-01  & -3.46e-01 \\ 
68 & 7 &  2 & 9.79e-01 &  2.91e-01  & 6.03e+00 &  1.09e-02 &  3.36e-01  & 5.39e-02 \\ 
69 & 9 &  4 & 4.48e-01 &  3.56e-01  & 7.49e+00 &  6.76e-01 &  1.91e+00  & 1.75e+00 \\ 
70 & 9 &  4 & 5.97e-01 &  6.35e-01  & 7.47e+00 &  1.75e-01 &  2.49e+00  & 1.80e+00 \\ 
71 & 3 &  4 & 1.12e+00 &  1.12e+00  & 6.34e+00 &  -3.81e+00 &  -1.56e+00  & -7.92e-01 \\ 
72 & 5 &  6 & 1.83e+00 &  1.01e+00  & 7.72e+00 &  -2.92e+00 &  1.25e-01  & 5.67e-01 \\ 
73 & 8 &  9 & 3.06e-01 &  5.78e-01  & 1.07e+01 &  1.33e+00 &  2.85e+00  & 1.25e+00 \\ 
74 & 16 &  15 & 1.01e+00 &  7.50e-01  & 1.58e+01 &  -9.38e-01 &  -8.88e-01  & -9.14e-01 \\ 
75 & 6 &  4 & 6.04e-01 &  4.13e-01  & 6.87e+00 &  2.04e-01 &  4.50e-01  & 6.32e-01 \\ 
76 & 6 &  4 & 6.07e-01 &  6.17e-01  & 6.96e+00 &  -1.63e+00 &  6.41e-02  & 3.97e+00 \\ 
77 & 6 &  4 & 2.24e-01 &  6.94e-01  & 7.01e+00 &  1.82e+00 &  -1.64e-01  & 2.60e+00 \\ 
78 & 3 &  4 & 8.93e-01 &  5.62e-01  & 6.26e+00 &  -2.29e+00 &  -4.09e+00  & -1.35e+00 \\ 
79 & 8 &  16 & 2.13e+00 &  2.08e+00  & 1.29e+01 &  -8.31e-03 &  -3.44e-02  & -3.32e-03 \\ 
80 & 3 &  2 & 5.73e-01 &  4.53e-01  & 5.35e+00 &  -3.26e-01 &  -4.45e-01  & -5.54e-01 \\ 
91 & 5 &  6 & 6.89e-01 &  5.08e-01  & 7.65e+00 &  -1.31e+00 &  -2.37e+00  & -7.94e-01 \\ 
92 & 5 &  6 & 8.00e-01 &  7.50e-01  & 7.20e+00 &  -8.88e+00 &  -2.92e+00  & -2.87e+00 \\ 
93 & 5 &  8 & 5.55e-01 &  8.38e-01  & 9.61e+00 &  -1.38e+00 &  -8.79e-01  & -1.25e+00 \\ 
94 & 5 &  8 & 8.42e-01 &  1.11e+00  & 1.26e+01 &  -6.12e+00 &  -7.87e+00  & -5.91e+00 \\ 
95 & 5 &  8 & 8.91e-01 &  1.21e+00  & 1.18e+01 &  -2.38e+00 &  -2.22e+00  & -1.43e+00 \\ 
96 & 5 &  8 & 6.57e-01 &  7.67e-01  & 9.44e+00 &  -5.38e+00 &  -4.99e+00  & -1.94e+00 \\ 
97 & 5 &  8 & 1.74e+00 &  6.35e-01  & 9.43e+00 &  -1.41e+00 &  -8.03e-01  & -7.87e-01 \\ 
98 & 5 &  8 & 5.12e-01 &  6.16e-01  & 9.51e+00 &  -3.33e+00 &  -9.37e+00  & -2.85e+00 \\ 
99 & 5 &  4 & 7.45e-01 &  9.34e-01  & 6.62e+00 &  -8.50e-02 &  -2.06e-01  & -1.19e-01 \\ 
100 & 20 &  2 & 1.32e+00 &  8.71e-01  & 5.95e-01 &  -8.90e-03 &  -9.19e-03  & -8.08e-03 \\ 
106 & 6 &  8 & 5.14e-01 &  7.57e-01  & 9.92e+00 &  -4.53e-02 &  -7.45e-02  & -5.13e-02 \\ 
107 & 5 &  3 & 1.71e+00 &  8.67e-01  & 4.14e+00 &  -8.33e-04 &  9.42e-06  & 6.46e-03 \\ 
108 & 10 &  4 & 1.63e+00 &  9.84e-01  & 7.81e+00 &  -5.63e-01 &  -6.02e-01  & -1.46e-01 \\ 
109 & 40 &  4 & 9.68e-01 &  1.94e+00  & 3.93e+01 &  -5.13e-03 &  -5.40e-03  & -4.98e-03 \\ 
110 & 40 &  4 & 1.17e+00 &  2.85e+00  & 1.76e+00 &  -5.12e-03 &  -5.53e-03  & -5.00e-03 \\ 
114 & 48 &  1 & 4.43e+00 &  1.25e+01  & 5.27e+01 &  -2.64e-01 &  -2.66e-01  & -2.84e-01 \\ 
115 & 9 &  4 & 6.30e-01 &  7.74e-01  & 7.62e+00 &  -1.18e-02 &  -8.45e-03  & 1.13e-01 \\ 
116 & 10 &  6 & 3.01e-01 &  2.03e-01  & 8.58e+00 &  4.80e-02 &  -1.43e-03  & 1.81e-02 \\ 
117 & 11 &  16 & 6.07e-01 &  1.08e+00  & 1.69e+01 &  5.21e-01 &  6.97e-01  & 8.40e-01 \\ 
118 & 9 &  4 & 7.13e-01 &  5.06e-01  & 7.58e+00 &  -1.52e-02 &  -3.33e-03  & 1.42e-03 \\ 
119 & 7 &  15 & 4.50e-01 &  1.04e+00  & 1.41e-01 &  3.64e-01 &  2.48e-01  & 3.61e-01 \\ 
120 & 5 &  9 & 4.04e-01 &  3.77e-01  & 1.00e+01 &  2.45e-01 &  3.36e-01  & 6.90e-01 \\ 
121 & 5 &  6 & 7.38e-01 &  6.77e-01  & 7.74e+00 &  -4.56e-02 &  -3.06e-01  & -8.28e-04 \\ 
122 & 9 &  16 & 1.74e+00 &  1.18e+00  & 1.49e+01 &  3.98e-04 &  -2.59e-03  & 5.18e-01 \\ 
123 & 6 &  9 & 4.71e-01 &  8.40e-01  & 1.03e+01 &  -3.13e-02 &  -2.05e-02  & 1.76e-02 \\ 
124 & 6 &  8 & 5.03e-01 &  5.10e-01  & 9.56e+00 &  1.73e+00 &  1.37e+00  & 7.70e-02 \\ 

\end{longtable}


\section{Nonlinear Eigenvalue Problems}\label{s:nonlinear}
In this section we consider time-delay systems of the form,
\[
v'(t) = \sum_{j=0}^m A_j(x)v(t-\tau_j).
\]

In this case, we have a nonlinear eigenvalue problem. To solve for the 
eigenvalues, we find the solutions $\lambda(x)$ of,
\[
\det (\Lambda(\lambda;x))=0,
\]
with,
\[
\Lambda (\lambda;x) = \lambda I -A_0(x)-\sum_{j=1}^m A_j(x) e^{-\lambda \tau_j}.
\]

The number of eigenvalues in this case is generally infinite, but within any right half-plane the number 
of eigenvalues is finite~\cite{Michiels2007}. In the numerical experiments, 
we find all of those which are to the right of $r=-1/\tau_m$, where $\tau_m$ is the maximum 
of the delays. If none are found, we repeatedly double $r$ until at least one eigenvalue 
appears. 

It can be shown that, in the case where each $A_i(x)$ depends smoothly on $x$ and where the eigenvalue has 
multiplicity 1, the derivative of the surface corresponding to each eigenvalue 
$\lambda_i$ is equal to~\cite{Michiels2007},
\[
\nabla_x \lambda_i = \frac{u^*_i \left(\frac{\partial A_0}{\partial x}+\sum_{j=1}^m \frac{\partial A_j}{\partial x}e^{-\lambda_i \tau_j}\right) v_i}{u^*_i \left(I+\sum_{j=1}^m\tau_j e^{-\lambda \tau_j} A_j\right) v_i}.
\]

Note that the term in the numerator, $\frac{\partial A_0}{\partial x}+\sum_{j=1}^m \frac{\partial A_j}{\partial x}e^{-\lambda_i \tau_j}$ corresponds to $\nabla_x F(x)$, and so we can use the same formula for the second derivative,~\eqref{eq:secondderiv} as in the linear case, with $F(x)$ replaced by $F(x,\lambda)=A_0(x)+\sum_{j=1}^m A_j(x) e^{-\lambda \tau_j}$. 

Alternatively, second-derivatives may be calculated explicitly. Specifically, 
\[
\begin{array}{l}
\nabla^2_{xx} \lambda_i(x) = -\frac{u_i^*(\nabla^2_{x\lambda}\Lambda(\lambda_i,x) \otimes \nabla_x \lambda_i+\nabla^2_{xx}\Lambda(\lambda_i,x) +\nabla^2_{\lambda\lambda}\Lambda(\lambda_i,x) \otimes (\nabla_x \lambda_i)
\otimes (\nabla_x \lambda_i)) v_i}{u_i^* \nabla_\lambda \Lambda(\lambda_i,x) v_i}
\\ \qquad +\frac{u_i^* (2\nabla_x \Lambda(\lambda_i,x)+2\nabla_\lambda\Lambda (\lambda_i,x)\otimes \nabla_x \lambda_i) \nabla_x v_i}{u_i^* \nabla_\lambda \Lambda(\lambda_i,x) v_i},
\end{array}
\]
where $\nabla_x v_i$ can be calculated (along with $\nabla_x \lambda_i$) by,
\[
\begin{pmatrix}
\Lambda(\lambda_i,x) & \nabla_\lambda \Lambda(\lambda_i,x) \\ 
2 v_i^* & 0
\end{pmatrix}
\begin{pmatrix}
\nabla_x v_i \\ \nabla_x \lambda_i 
\end{pmatrix}
=\begin{pmatrix}
\nabla_x v_i \\ 0
\end{pmatrix},
\]
where the second set of equations comes from differentiating $v_i^* v_i = 1$.

We compare the performance of 500 trials of SLP and HANSO with and without gradient sampling for 
two time-delay systems described in~\cite{Van2008}.

The first example is a third-order feedback controller system of the form,
\[
v'(t) = A v(t)+B(x) v(t-5),
\]
with $A$ and $B(x)$ defined to be,
\[
A = \begin{pmatrix}

-0.08 & -0.03 & 0.2 \\ 0.2 & -0.04 & -0.005 \\ -0.06 & -0.2 & -0.07 
\end{pmatrix}
\]
and 
\[
B(x) = \begin{pmatrix} -0.1 \\ -0.2 \\ 0.1 \end{pmatrix}
\begin{pmatrix} x_1 & x_2 & x_3 \end{pmatrix}.
\]

We present the results below in Table~\ref{tablenonlin1}. In this case, statistically, SLP 
finds no lower or higher minimizer than HANSO with or without gradient sampling. However, 
it almost always finds it in less time. In this case, SLP appears to also perform similarly in terms 
of finding the lowest minimizer as SQP, but now takes less time. We also include the mean and standard 
deviation of 
the values and times in Table~\ref{tablenonlin1b}. One can consider the true mean time of execution and 
value at the final solution as a distribution dependent on the initial starting point with unknown form 
for each algorithm. Therefore, we can take a sample difference of means test to compare the performance. 
If we take a Student t-test for the difference in means for the times, then even for $0.0001$ significance, the difference 
in both mean times and values is significant between SLP and HANSO (both with and without gradient sampling). 

\begin{table}
\caption{\label{tablenonlin1} Number of times SLP outperformed HANSO and SQP (out of 500 sample runs). 
Value for each solver for each run is taken as the best of 10 random starting points. 
Time is the total clock time taken to perform the ten runs.}
\begin{center}
\begin{tabular}{l  c  r } 
\toprule 
 & in value  & in time \\ \midrule 
HANSO & 250 & 500   \\ 
HANSO without gradient sampling & 257 & 447  \\
SQP & 227 & 344 \\
\bottomrule
\end{tabular}
\end{center}
\end{table}

\begin{table}
\caption{\label{tablenonlin1b} Mean (standard deviation) for values and times for HANSO, HANSO without gradient 
sampling, SLP, and SQP (out of 500 sample runs). Value for each solver for each run is taken as the best of 10 random starting points. 
Time is the total clock time taken to perform the ten runs.}
\begin{center}
\begin{tabular}{l  c  r } 
\toprule 
 & value  & time \\ \midrule 
HANSO & -0.074 (0.036) & 177 (30.3)   \\ 
HANSO without gradient sampling & -0.069 (0.036) & 6.1 (2.7)  \\
SLP &-0.081 (0.053) & 4.6 (1.1) \\
SQP & -0.088 (0.055) & 5.3 (1.6) \\
\bottomrule
\end{tabular}
\end{center}
\end{table}

The next example is given below,
\[
\begin{array}{l}
T_h \dot{x}_h(t) = -x_h(t-\eta_h)+K_b x_a(t-\tau_b)+K_u x_{h,set}(t-\tau_u), \\
T_a \dot{x}_a(t) = -x_a(t)+x_c(t-\tau_e)+K_a(x_h(t)-\frac{1+q}{2}x_a(t)-\frac{1-q}{2}x_c(t-\tau_e)), \\
T_d \dot{x}_d(t) = -x_d(t)+K_d x_a(t-\tau_d), \\
T_c \dot{x}_c(t)  = -x_c(t-\eta_c)+K_c x_d(t-\tau_c), \\
\dot{x}_e(t) = -x_c(t)+x_{c,set}(t),
\end{array}
\]
with,
\[
x_{h,set}(t) = \begin{pmatrix} K_1 & K_2 & K_3 & K_4 & K_5 \end{pmatrix} 
\begin{pmatrix} x_h(t) & x_a(t) & x_d(t) & x_c(t) & x_e(t) \end{pmatrix}^T. 
\]

The results, which are qualitatively similar as in the first example, are given in Tables~\ref{tablenonlin2}
and~\ref{tablenonlin2b}.

\begin{table}
\caption{\label{tablenonlin2} Number of times SLP outperformed HANSO (out of 500 sample runs).}
\begin{center}
\begin{tabular}{l  c  r } 
\toprule 
 & in value  & in time \\ \midrule 
HANSO & 268 & 500   \\ 
HANSO without gradient sampling & 283 & 466  \\
\bottomrule
\end{tabular}
\end{center}
\end{table}

\begin{table}
\caption{\label{tablenonlin2b} Mean (standard deviation) for values and times for HANSO, HANSO without gradient 
sampling, SLP, and SQP (out of 500 sample runs).}
\begin{center}
\begin{tabular}{l  c  r } 
\toprule 
 & value  & time \\ \midrule 
HANSO & -0.077 (0.0052) & 2,200 (6400)   \\ 
HANSO without gradient sampling & -0.067 (0.0054) & 69.0 (28)  \\
SLP &-0.083 (0.0062) & 83.5 (140) \\
SQP & -0.088 (0.0103) & 75.5 (76) \\
\bottomrule
\end{tabular}
\end{center}
\end{table}

\section{Conclusion}\label{s:conclusion}
In this paper we studied the eigenvalue optimization problem of minimizing the spectral abscissa, the maximum 
real eigenvalue part. This problem is important for designing stabilizing controllers, such as in LTI and time-delay models.
We presented 
an algorithm that incorporated linear and quadratic models of eigenvalue surfaces corresponding to different 
eigenvalues in a sequential linear and sequential quadratic programming framework. We expected this 
to produce a faster and possibly more reliable algorithm for finding minima of the spectral abscissa, since 
the model is capable of approximating the spectral abscissa surface past points of nonsmoothness. 

Our numerical results for both linear and nonlinear problems borne out our expectation, and we find that, in particular, 
the sequential linear variant of the algorithm tends to, in general, outperform the most comparable competitor, HANSO. 
For linear problems, it is faster and more reliable than HANSO with gradient sampling, and more reliable but slower 
than HANSO 
without gradient sampling. For nonlinear problems, SLP is faster than and at least equally as reliable as both
HANSO and HANSO without gradient sampling. 

In future research we intend to study the pseudospectral abscissa, which is a locally Lipschitz function, and hence 
more readily amenable to theoretical convergence analysis. In addition, researh can include
testing these algorithms on interesting large-scale delay and control applications.

\section{Acknowledgements}
This research was supported by Research Council KUL: 
PFV/10/002 Optimization in Engineering Center OPTEC, GOA/10/09 MaNet,
 Belgian Federal Science Policy Office: IUAP P7 (DYSCO, Dynamical systems, 
control and optimization, 2012-2017); ERC ST HIGHWIND (259 166).
V. Kungurtsev was also supported by
        the European social fund within the framework of realizing the project ``Support of inter-sectoral
        mobility and quality enhancement of research teams at the Czech  Technical University in Prague'',
        CZ.1.07/2.3.00/30.0034.

\bibliographystyle{myplain}
\bibliography{eigenval}

\end{document}